\documentclass[a4paper]{article}
\usepackage{amsmath,amsfonts,amssymb}
\usepackage{multicol,multirow,diagbox,bm}
\usepackage{graphicx,subfigure}
\usepackage{lineno,cases}
\allowdisplaybreaks 
\usepackage{pgf,tikz}
\usetikzlibrary{arrows,shapes,automata,positioning}
\usetikzlibrary{fit}
\usepackage{pst-node}
\usepackage{cite}
\usepackage{color,xcolor}
\renewcommand{\bar}{\overline}
\renewcommand{\top}{\intercal}

\newcommand{\e}{\varepsilon}
\newcommand{\R}{\mathbb{R}}

\usepackage{bm} 
\newcommand{\N}{\mathcal{N}}
\newcommand{\G}{\mathcal{G}}
\newcommand{\E}{\mathcal{E}}
\newcommand{\A}{\mathcal{A}}

\newtheorem{ass}{Assumption}
\newtheorem{lemma}{Lemma}

\newtheorem{corollary}{Corollary}
\newtheorem{defn}{Definition}
\newtheorem{theorem}{Theorem}
\newtheorem{remark}{Remark}

\newcommand{\pb}{\noindent\textbf{Proof.} }
\newcommand{\pe}{\hfill\rule{4pt}{8pt}}

\usepackage{cite}
\begin{document}
\title{ Primal-dual $\varepsilon$-Subgradient Method for Distributed Optimization 
	\thanks{This work was supported by National Natural Science Foundation of China under  Grants 61973043.}
}
\author{Kui Zhu and Yutao Tang
	\thanks{K. Zhu and Y. Tang are both with the School of Artificial Intelligence, Beijing University of Posts and Telecommunications, Beijing 100876, China (e-mails: kuizhu\_19@butp.edu.cn, yttang@bupt.edu.cn). }
}	
\date{ }
\maketitle

{\noindent\bf Abstract}:  This paper studies the distributed optimization problem when the objective functions might be nondifferentiable and subject to heterogeneous set constraints.  Unlike existing subgradient methods, we focus on the case when the exact subgradients of the local objective functions can not be accessed by the agents. To solve this problem, we propose a projected primal-dual dynamics using only the objective function's approximate subgradients. We first prove that the formulated optimization problem can generally be solved with an error depending upon the accuracy of the available subgradients.  Then, we show the exact solvability of this distributed optimization problem when the accumulated approximation error of inexact subgradients is not too large. After that, we also give a novel componentwise normalized variant to improve the transient behavior of the convergent sequence.  The effectiveness of our algorithms is verified by a numerical example.
	

\section{Introduction}

The last decade has witnessed considerable interests in distributed optimization problems due to the numerous applications in signal processing, control, and machine learning. To solve this problem,  subgradient information of the objective functions has been widely used due to the cheap iteration cost and well-established convergence properties \cite{nedic2018distributed,yang2017distributed}. 

Note that most of these subgradient-based results assume the availability of local cost function's exact subgradients. In many circumstances, the function subgradient is computed by solving another auxiliary optimization problem as shown in \cite{bertsekas2015convex,kiwiel2004convergence,devolder2014first}. In practice, we are often only able to solve these subproblems approximately. Hence, in that context, numerical methods solving the original optimization problem are provided with only inexact subgradient information. This leads us to investigate the solvability of distributed optimization problem using inexact subgradient information.

A close topic is the inexact augmented Lagrangian method. As surveyed in \cite{jakovetic2020primal}, this method has been extensively extended to distributed settings in various ways assuming the primal variables are only obtained in an approximate sense.  Nevertheless, most of these results  still require the exact gradient or subgradient information of the local objective functions at each given estimate. It is interesting to ask whether the primal-dual method is still effective when only inexact gradient/subgradient information is available.

In this paper, we focus on a typical distributed consensus optimization problem for a sum of convex objective functions subject to heterogeneous set constraints. Although this problem has been partially studied by gradient/subgradient methods in \cite{nedic2009distributed, jakovetic2014linear,yi2015distributed, lei2016primal, xi2016distributed, liu2017convergence, zeng2017distributed, zhu2020boundedness}, its solvability using only inexact subgradient information has not yet been addressed and is still unclear at present.

To solve this problem,  we first convert it into a distributed saddle point seeking problem and present a projected primal-dual $\e$-subgradient dynamics to handle both the distributedness and constraints. When the objective functions are smooth with exact gradients, the proposed algorithms reduce to the primal-dual dynamics considering in \cite{lei2016primal,jakovetic2014linear}. Then, we discuss the convergent property  with diminishing step sizes and suboptimality of the proposed algorithm depending on the accuracy of available subgradients. In particular, we show that if the accumulated error resulting from the subgradient inexactness is not too large, the proposed algorithm under certain diminishing step size will drive all the estimates of agents to reach a consensus about an optimal solution to the global optimization problem.  To our knowledge, this might be the first attempt to solve the formulated distributed optimization problem using only inexact subgradients of the local objective functions.   To improve the transient performance of our preceding designs, we further propose a novel componentwise normalized step size as that in \cite{zhu2020boundedness}.  As a byproduct, this normalized step size  removes the widely used subgradient boundedness assumption in the literature \cite{nedic2009distributed, jakovetic2014linear}.  


The rest of this paper is organized as follows. We first give some preliminaries in Section \ref{sec:pre} and then introduce the formulation of our problem in Section \ref{sec:form}. Main results are presented in Section  \ref{sec:main}. After that, we give a numerical example in Section \ref{sec:simu} to show the effectiveness of our design.  Finally, some concluding remarks are given in Section \ref{sec:con}.

\section{Preliminary}\label{sec:pre}
In this section, we first give some preliminaries about graph theory and convex analysis.

\subsection{Graph Theory}
Let $\R^n$ be the $n$-dimensional Euclidean space and $\R^{n\times m}$ be the set of all $n\times m$ matrices. ${\bf 1}_n$ (or ${\bf 0}_n$) denotes an $n$-dimensional all-one (or all-zero) column vector and ${\bm 1}_{n\times m}$ (or ${\bm 0}_{n\times m}$) all-one (or all-zero) matrix.  $\mathrm{col}(a_1,\,{\dots},\,a_n) = {[a_1^\top,\,{\dots},\,a_n^\top]}^\top$ for column vectors $a_1,\,\dots, a_n$.  For a vector $x$ (or matrix $A$) , $||x||$ ($||A||$) denotes its Euclidean (or spectral) norm. 

A weighted (undirected) graph $\mathcal {G}=(\mathcal {N}, \mathcal {E}, \mathcal{A})$ is defined as follows, where $\mathcal{N}=\{1,{\dots},n\}$ is the set of nodes, $\mathcal {E}\subset \mathcal{N}\times \mathcal{N}$ is the set of edges, and $\mathcal{A}\in \mathbb{R}^{n\times n}$ is a weighted adjacency matrix. $(i,j)\in \mathcal{E}$ denotes an edge leaving from node $i$ and entering node $j$. The weighted adjacency matrix of this graph $\mathcal {G}$ is described by $\mathcal{A}=[a_{ij}]\in \R^{n\times n}$, where $a_{ii}=0$ and $a_{ij}\geq 0$ ($a_{ij}=a_{ji}>0$ if and only if there is an edge between nodes $j$ and $i$\,).  The neighbor set of agent $i$ is defined as $\mathcal{N}_i=\{j\colon (j,\,i)\in \mathcal {E} \}$ for $i=1,\,\dots,\,n$. A path in graph $\mathcal {G}$ is an alternating sequence $i_{1}e_{1}i_{2}e_{2}{\cdots}e_{k-1}i_{k}$ of nodes $i_{l}$ and edges $e_{m}=(i_{m},\,i_{m+1})\in\mathcal {E}$ for $l=1,\,2,\,{\dots},\,k$. If there exists a path from node $i$ to node $j$ then node $i$ is said to be reachable from node $j$. 
The Laplacian $L=[l_{ij}]\in \mathbb{R}^{n\times n}$ of graph $\mathcal{G}$ is defined as $l_{ii}=\sum_{j\neq i}a_{ij}$ and $l_{ij}=-a_{ij} (j\neq i)$. It can be found that the Laplacian  is symmetric and semi-definite. Denote its ordered eigenvalues as $0=\lambda_1\leq \lambda_2\leq \dots\lambda_N$. The corresponding eigenvector of $\lambda_1=0$ is the all one vector ${\bm 1}_N$. Moreover, $\lambda_2>0$ if and only if this graph $\mathcal{G}$ is connected.

\subsection{Convex analysis}

For each set $X\in \R^m$, the indicator function is denoted by $\delta_{X}$ with $\delta_{X}(x)=0$ for any $x\in X$ and $\delta_{X}(x)=\infty$ for any $x\notin X$.  A set $X\in \R^m$ is said to be convex if $\theta x+(1-\theta)y \in X$ for any any $x,\,y\in X$ and $\theta\in (0,\,1)$. For a closed convex set $X\neq \emptyset$, the projection operator $P_{X}\colon \R^m\to X$ is defined as $P_{X}[x]=\arg\min_{y\in X}||y-x||$. The projection operator is non-expansive in the sense that $||P_{X}[x]-P_{X}[y]||\leq ||x-y||$ for any $x,\,y\in \R^m$. For any $x\in \R^m$, it holds $(x-P_{X}[x])^\top (y-P_{X}[x])\leq 0,\,\forall y \in X$.

For a given function $f\colon \R^m\to \R$,  we denote by $\mathrm{dom}\,f=\{x\in \R^m \mid |f(x)|<\infty \}$ the domain of $f$. We always assume that $\mathrm{dom}\,f\neq \emptyset$ and $\mathrm{dom}\,f=\R^m$ if not specified. We say it is convex if its domain is convex and  $f(\alpha x+(1-\alpha y))\leq \alpha f(x)+(1-\alpha)f(y)$ holds for all $x,\,y\in \mathrm{dom}\,f$ and $\alpha\in [0,\,1]$.  If this inequality is strict in the sense that the equation holds only if $x=y$, the function is called strictly convex.  A function $f$ is called closed and convex on a convex set $X\subset \mathrm{dom}\,f$ if its constrained epigraph $\mbox{epi}_{X}(f)=\{(x,\,t)\in X\times \R\mid t\geq f(x) \}$ is a closed convex set. If $X=\mathrm{dom}\,f$, we call $f$ a closed convex function. 


A vector-valued function ${\bm f}\colon \R^m \rightarrow \R^m$ is Lipschitz with constant $\vartheta>0$ (or simply $\vartheta$-Lipschitz) if we have
\begin{align*}
	||{\bm f}(\zeta_1)-{\bm f}(\zeta_2)||\leq \vartheta ||\zeta_1-\zeta_2||,\, \forall \zeta_1, \zeta_2 \in \R^m
\end{align*}

Let us consider a function $\phi\colon X\times Z\to \R$, where $X$ and $Z$ are nonempty subsets of $\R^n$ and $\R^m$, respectively. A pair of vectors $x^*\in X$ and $z^*\in Z$ is called a saddle point of $\phi$ if $\phi(x^*,\,z)\leq \phi(x^*,\,z^*)\leq \phi(x,\,z^*)$ holds for any $x\in X$ and $z\in Z$. 

\section{Problem Formulation}\label{sec:form}

In this paper, we focus on solving the following constrained optimization problem by a network of $N$ agents: 
\begin{align}\label{prop:main}
	\begin{split}
		\min&\quad f(x)=\sum_{i=1}^N f_i(x) \\
		\mbox{s.t.}& \quad x  \in X\triangleq \bigcap_{i=1}^N X_i
	\end{split}
\end{align}
Here function $f_i\colon \R\to \R$ and set $X_i$ are private to agent $i$ for each $i\in \mathcal{N}\triangleq \{1,\,2,\, \dots,\,N\}$ and can not be shared  with others. 

To ensure its solvability, the following assumption is made. 
\begin{ass}\label{ass:convex}
	For each $i\in \N$,  function $f_i\colon \R\to \R$ is  convex, set $X_i$ is convex and closed, and $\mathrm{int}\,X=\bigcap_{i\in \N}\mathrm{int}\,X_i$ is nonempty and contained in $\mathrm{dom}\,f_i$. 
\end{ass}

Note that $f_i$ might be nondifferentiable under this assumption.  Denote the minimal value of  problem \eqref{prop:main} by $f^{*}$ and the optimal solution set by $\cal{X}^*$, i.e., $f^*=\min_{x\in X} f(x)$ and $\mathcal{X}^*=\{x \in X \mid f(x)=f^*\}$.  As usual, we assume $f^*$ is finite and set $\mathcal{X}^*$ is nonempty.  To cooperatively address the  optimization problem \eqref{prop:main} in a distributed manner, we use a weighted undirected graph $\mathcal{G}=\{\N,\,\,\E,\, \A\}$ to describe the information sharing relationships with node set $\N$, edge set $\E\subset \N\times \N$, and weight matrix $\A=[a_{ij}]_{N\times N}$.  Here $a_{ij}=a_{ji}>0$ means agents $i$ and $j$ can communicate with each other.

\begin{ass}\label{ass:graph}
	Graph $\mathcal{G}$ is connected.
\end{ass}

Suppose agent $i$ maintains an estimate $x_i$ of the optimal solution to  \eqref{prop:main} with  other (possible) auxiliary variables. Agents exchange these variables through the communication network described by $\G$ and perform some updates at given discrete-time instants $k= 1,\,2,\,\dots$.  
Then, the distributed optimization problem in this paper is formulated to find an update rule of $x_i(k)$  for agent $i$ using only its own and neighboring information such that $\lim_{k\to \infty}[x_i(k)-x_j(k)]=0$ for any $i,\,j\in \mathcal{N}$ and $\lim_{k\to\infty}\sum_{i=1}^N f_i(x_i(k))=f^*$. If possible, we expect that all the estimates will converge to an optimal solution to problem \eqref{prop:main}.

As stated above, this problem has been intensively studied in literature \cite{nedic2018distributed,yang2017distributed}. However, most existing designs require the exact subgradients of the local objective functions in constructing effective distributed algorithms. In this paper, we are interested in the solvability of the formulated  distributed constrained optimization problem \eqref{prop:main} working with inexact subgradients of the local objective functions.  For this purpose, we adopt the notion of $\e$-subgradient to describe such inexactness as in \cite{polyak1987introduction} and assume the $\varepsilon$-subgradient of $f_i$ can be easily computed for any given $\e\geq 0$. 

\begin{defn}
	For a convex function $f\colon \R^m \to \R$ and a scalar $\varepsilon>0$,    $g\in \R^m $ is said to be an $\varepsilon$-subgradient of $f$ at $x\in \R^m$ if  $$f(y) \geq f(x)+g^\top (y-x)-\varepsilon,\quad \forall y \in \mathrm{dom}\,f$$ 
\end{defn}

Denote by $\partial _{\varepsilon} f(x)$ the $\varepsilon$-subdifferential of $f$ at $x\in \R^m$, which is the set of all $\varepsilon$-subgradients of $f$ at $x$.  $\partial _{\varepsilon} f(x)$ is nonempty and convex for any $x\in \R^m$ due to the convexity of $f$. Moreover, $\partial_{0} f(x)$ coincides with the subdifferential of $f$ at $x\in \R^m$.

In next section, we will convert our  problem into a saddle-point seeking problem and develop a projected primal-dual $\varepsilon$-subgradient method with rigorous solvability analysis.

\section{Main Result}\label{sec:main}

To begin with, we rewrite problem \eqref{prop:main} into an alternative form as that in \cite{lei2016primal,liu2017constrained}: 
\begin{align}\label{prop:main-reform}
	\min~~&\tilde f({\bf x})=\sum_{i=1}^N f_i(x_i) \nonumber \\
	\mbox{s.t. }\, & L{\bf x}={\bf 0}_N\\
	&{\bf x}\in \tilde X\triangleq X_1\times \dots\times X_N \nonumber
\end{align}
where ${\bf x}=\mathrm{col}(x_1,\,\dots,\,x_N)$ and $L$ is the Laplacian of graph $\G$. Note that $L$  is symmetric and positive semi-definite with  its ordered eigenvalues as $0=\lambda_1<\lambda_2\leq \dots\lambda_N$  under Assumption \ref{ass:graph} by Theorem 2.8 in \cite{mesbahi2010graph}.

Consider the augmented Lagrangian function of problem \eqref{prop:main-reform}:
\begin{align}\label{func:lagrangian}
	\Phi({\bf x},\, {\bf v})=\tilde f({\bf x})+{\bf v}^\top L {\bf x}+\frac{1}{2}{\bf x}^\top L{\bf x}
\end{align}
with ${\bf v} =\mbox{col}(v_1,\,\dots,\,v_N)\in \R^N$.  By Proposition 3.4.1 in \cite{bertsekas2015convex}, if  $\Phi$ has a saddle point $({\bf x}^*,\,{\bf v}^*)$ in $\tilde X\times \R^N$, then ${\bf x}^*$ must be an optimal solution to problem \eqref{prop:main-reform}, which in turn provides an optimal solution to \eqref{prop:main}. 
Since the Slater's condition holds under Assumption \ref{ass:convex}, such saddle points indeed exist by virtue of Theorems 3.34 and 4.7 in \cite{ruszczynski2011nonlinear}. Thus, it suffices for us to seek a saddle point  of $\Phi$ in $\tilde X\times \R^N$. 

Following this conversion, many solvability results on problem \eqref{prop:main}  have been presented when the exact gradient or subgradient of $f_i$ is available, e.g., \cite{wang2010control, gharesifard2014distributed, liu2017constrained, lei2016primal,kia2015distributed, zhu2020boundedness}.  However, whether and how $\varepsilon$-subgradient algorithms can be derived  has not been discussed yet. To this end, we are motivated by aforementioned saddle-point seeking designs and present the following dynamics:
\begin{align}\label{alg:constrained}
	\begin{split}
		x_{i}(k+1) &= P_{X_i}[x_{i}(k) - \alpha_k (g_{i}(k)+ \hat x_i(k)+\hat v_i(k))] \\
		v_{i}(k+1) &=v_{i}(k)+ \alpha_k \hat x_i(k)
	\end{split}
\end{align}
where $\hat x_i(k) \triangleq \sum_{j=1}^N a_{ij}(x_i(k)-x_j(k))$, $\hat v_i(k) \triangleq \sum_{j=1}^N a_{ij}(v_i(k)-v_j(k))$,  and $g_i(k)\in \partial_{\e_k} f_i(x_i(k))$ with parameters $\e_k,\,\alpha_k>0$ to be specified later. 
It can be taken as a constrained version of algorithms in \cite{wang2010control, gharesifard2014distributed}. 
Different from  similar primal-dual designs in \cite{lei2016primal,liu2017constrained}, we do not require the differentiability  of these objective functions or their exact gradients. 

Letting ${\bf x}(k)=\mathrm{col}(x_1(k), \dots, x_N(k))$ and ${\bf v}(k)=\mathrm{col}(v_1(k), \dots, v_N(k))$, we can put \eqref{alg:constrained} into a compact form:
\begin{align}\label{alg:constrained-compact}
	\begin{split}
		{\bf x}(k+1)&=P_{\tilde X}[{\bf x}(k)- \alpha_k ({\bf g}(k)+L{\bf v}(k)+L{\bf x}(k))]\\
		{\bf v}(k+1)&={\bf v}(k)+ \alpha_k L{\bf x}(k)
	\end{split}
\end{align}
with 
${\bf g}(k)=\mathrm{col}(g_1(k),\,\dots,\,g_N(k))\in \partial_{N\e_k}\tilde f({\bf x}(k))\in \R^N$. It can be further rewritten as follows.
\begin{align}\label{alg:constrained-compact-proof}
	{\bf z}(k+1)=P_{\bar X}[{\bf z}(k)-\alpha_kT_{\e_k}({\bf z}(k))]
\end{align}
where  ${\bf z}(k)=\mathrm{col}({\bf x}(k),\,{\bf v}(k))$, $\bar X=\tilde X\times \R^N$, and 
\begin{align*}
	T_{\e_k}({\bf z}(k))=\begin{bmatrix}
		{\bf g}(k)+L{\bf v}(k)+L{\bf x}(k)\\
		-  L{\bf x}(k)
	\end{bmatrix}
\end{align*} 

To establish the effectiveness of this algorithm, another assumption is made as follows.  
\begin{ass}\label{ass:boundedness}
	The $\e_k$-subgradient sequence $\{g_i(k)\}$ is uniformly bounded for each $i$, i.e., there exists a scalar $C>0$ such that $\max_{i\in \N}\{||g_i(k)||\}<C$ all $k>0$.
\end{ass}

This assumption is temporally made for simplicity as in \cite{nedic2009subgradient,jakovetic2014linear} and will be further removed later by some novel step sizes later. 
Suppose  ${\bf z}^*=\mathrm{col}({\bf x}^*, \,{\bf v}^*)$ is a saddle point of $\Phi$ in $\tilde X\times \R^N$. 
Here is a key lemma under Assumption \ref{ass:boundedness}. 
\begin{lemma}\label{lem:main}
	Suppose Assumptions \ref{ass:convex}--\ref{ass:boundedness} hold. Along the trajectory of algorithm \eqref{alg:constrained}, there exists some $C_1>0$ such that, for any $k=1,\,2,\,\dots$, the following inequality holds:
	\begin{align}\label{eq:lem:main}
		\begin{split}
			||{\bf z}(k+1)-{\bf z}^*||^2 &\leq (1+C_1\alpha_k^2)||{\bf z}(k)-{\bf z}^*||^2- 2\alpha_k\Delta({\bf x}(k)) +2N\alpha_k\e_k+C_1\alpha_{k}^2
		\end{split}
	\end{align}
	where $\Delta({\bf x}(k))\triangleq \Phi({\bf x}(k),\,{\bf v}^*)-\Phi({\bf x}^*,\,{\bf v}^*)+ \frac{1}{2}{\bf x}(k)^\top L{\bf x}(k)\geq 0$.
\end{lemma}
\pb  By lemma conditions, $({\bf x}^*, \,{\bf v}^*)$ is a saddle point of $\Phi$.  Then, $\Delta({\bf x}(k))\triangleq \Phi({\bf x}(k),\,{\bf v}^*)-\Phi({\bf x}^*,\,{\bf v}^*)+ \frac{1}{2}{\bf x}(k)^\top L{\bf x}(k)\geq 0$ can be easily verified by the definition of saddle points. 

Next, we consider the evolution of $||{\bf z}(k)-{\bf z}^*||^2$ with respect to $k$. Under the iteration \eqref{alg:constrained}, it follows then
\begin{align}\label{eq:thm:proof-error-eq1}
	&||{\bf z}(k+1)-{\bf z}^*||^2\nonumber\\
	&\qquad=||P_{\bar X}[{\bf z}(k)-\alpha_kT_{\e_k}({\bf z}(k))]-{\bf z}^*||^2 \nonumber \\
	&\qquad\leq ||{\bf z}(k)-\alpha_k T_{\e_k}({\bf z}(k))-{\bf z}^*||^2 -||{\bf z}(k)-\alpha_kT_{\e_k}({\bf z}(k))-P_{\bar X}[{\bf z}(k)-\alpha_kT_{\e_k}({\bf z}(k))]||^2 \nonumber \\
	&\qquad\leq ||{\bf z}(k)-{\bf z}^*||^2-2\alpha_k({\bf z}(k)-{\bf z}^*)^\top T_{\e_k}({\bf z}(k))+\alpha_{k}^2||T_{\e_k}({\bf z}(k))||^2 
\end{align}
By the proprieties of saddle point and $\e_k$-subgradient, we have
\begin{align}
	({\bf z}(k)-{\bf z}^*)^\top  T_{\e_k}({\bf z}(k))&= ({\bf x}(k)-{\bf x}^*)^\top ({\bf g}(k)+L{\bf v}(k)+L{\bf x}(k))- ({\bf v}(k)-{\bf v}^*)^\top L{\bf x}(k) \nonumber \\
	&\geq \tilde f({\bf x}^k)-\tilde f({\bf x}^*)-N\e_k+{{\bf v}^*}^\top L {\bf x}(k)+{\bf x}(k)^\top L{\bf x}(k) \nonumber\\
	&=\Phi({\bf x}(k),\,{\bf v}^*)-\Phi({\bf x}^*,\,{\bf v}^*)+ \frac{1}{2}{\bf x}(k)^\top L{\bf x}(k)-N\e_k  \label{eq:thm-proof-error-eq2}\\
	&=\Delta({\bf x}(k))-N\e_k \nonumber
\end{align}

Since $L {\bf x}^*={\bf 0}$, $T_{\e_k}({\bf z}(k))$ can be rewritten as 
\begin{align*}
	T_{\e_k}({\bf z}(k))=\begin{bmatrix}
		{\bf g}(k)+L({\bf x}(k)-{\bf x}^*)+L({\bf v}(k)-{\bf v}^*)+L{\bf v}^*\\
		- L({\bf x}(k)-{\bf x}^*)
	\end{bmatrix}
\end{align*} 
Under Assumption \ref{ass:boundedness}, there must be constant $C_1>0$ such that 
\begin{align}\label{eq:thm-proof-error-eq3}
	||	T_{\e_k}({\bf z}(k))||^2\leq C_1(1+||{\bf z}(k)-{\bf z}^*||^2)
\end{align}
Putting all inequalities \eqref{eq:thm:proof-error-eq1}--\eqref{eq:thm-proof-error-eq3} together, we have 
\begin{align*}
	||{\bf z}(k+1)-{\bf z}^*||^2 &\leq (1+C_1\alpha_k^2)||{\bf z}(k)-{\bf z}^*||^2 - 2\alpha_k\Delta({\bf x}(k)) +2N\alpha_k\e_k+C_1\alpha_{k}^2
\end{align*}
which is exactly the expected inequality \eqref{eq:lem:main}. 
\pe

When the exact subgradient is available (i.e., $\e_k=0$), the inequality \eqref{eq:lem:main} can be simplified into the well-known supermartingale inequality ensuring the convergence of ${\bf z}(k)$ towards ${\bf z}^*$ as shown in \cite{bertsekas2015convex} if $\{\alpha_k\}$ is chosen to satisfy
\begin{align}\label{eq:thm-parameter-condition-1}
	\sum_{k=1}^{\infty}\alpha_k=\infty,\quad \sum_{k=1}^{\infty}\alpha_k^2<\infty
\end{align}
However,  the inexactness of available subgradients deteriorates this property and the expected convergence might fail when we use only $\varepsilon$-subgradients in the iteration \eqref{alg:constrained}. 

Let us denote 
$\bar \Delta =  \liminf_{k\to \infty} \Delta({\bf x}(k))$  
and take a closer look at the inequality \eqref{eq:lem:main}. Note that $\bar \Delta $ consists of two parts, i.e., the discrepancy of function values and violation of constraints (in term of consensus error) of the iterative sequence. It can be taken as a measure of suboptimality of the iterative sequence. In other words, we can determine the upper bound for $\bar \Delta$ to evaluate the effectiveness of algorithm \eqref{alg:constrained}.  

We first consider the case when $\varepsilon_k$ is a constant.

\begin{theorem}\label{thm:main-delta}
	Suppose Assumptions \ref{ass:convex}--\ref{ass:boundedness} hold. Let the step size $\alpha_k$ be chosen to satisfy \eqref{eq:thm-parameter-condition-1} and $\varepsilon_{k}$ be fixed at some scalar $\varepsilon_0>0$. 
	Then, along the trajectory of algorithm \eqref{alg:constrained}, it holds that
	\begin{align}
		0\leq \bar \Delta\leq N\varepsilon_0
	\end{align}
\end{theorem}
\pb  To prove this theorem, we only have to show that $\bar \Delta\leq N\varepsilon_0$. If the inequality does not hold, there must exist a $\delta>0$ and a sufficient large integer $K_1>1$ such that
$
\Delta({\bf x}(k))> N\varepsilon_0+\delta
$ 
for all $k\geq K_1$. By \eqref{eq:thm-parameter-condition-1}, we have $\lim_{k\to \infty} \alpha_{k}=0$. Thus,  there must exist an integer $K_2>1$ such that  $0<\alpha_{k}\leq \frac{\delta}{C_1}$ for all $k\geq K_2$. Bringing these conditions together, one can strengthen inequality \eqref{eq:lem:main} for all $k\geq  K\triangleq \max\{K_1,\,K_2\}$ as follows.
\begin{align*}
	||{\bf z}(k+1)-{\bf z}^*||^2 &\leq (1+C_1\alpha_k^2)||{\bf z}(k)-{\bf z}^*||^2 -\alpha_k\delta
\end{align*}
Summing up its both sides from $K$ to $\bar K>K$ gives
\begin{align*}
	||{\bf z}({\bar K}+1)-{\bf z}^*||^2 & \leq ||{\bf z}(K)-{\bf z}^*||^2 \prod_{k=K}^{\bar K}(1+C_1\alpha_k^2)  - \delta \sum_{k=K}^{\bar K}\alpha_k
\end{align*}
where we use $1+C_1\alpha_k^2>1$ to handle the cross terms. 

Note that $1+\theta\leq e^\theta$ for any $\theta>0$. Hence $\prod_{k=K}^{\bar K}(1+C_1\alpha_k^2) \leq e^{C_1\sum_{k= K}^{\bar K}\alpha_k^2}\leq e^{C_1\sum_{k= 1}^{\infty}\alpha_k^2}$. Under the condition \eqref{eq:thm-parameter-condition-1}, there must exist a positive scalar $\bar C>0$ such that 
\begin{align*}
	||{\bf z}({\bar K}+1)-{\bf z}^*||^2 & \leq \bar C||{\bf z}(K)-{\bf z}^*||^2 - \delta \sum_{k=K}^{\bar K}\alpha_k
\end{align*}
which can not hold for a sufficiently large $\bar K$ since  $\sum_{k=1}^{\infty}\alpha_k=\infty$. We obtain a contradiction and complete the proof.
\pe

\begin{remark}
	According to Theorem \ref{thm:main-delta}, one can generally obtain a suboptimal solution to problem \eqref{prop:main} using inexact subgradients. If we are interested in an exact solution, it is required to ensure $\lim_{k\to \infty} \varepsilon_k=0$. For a very special case when $\e_k=0$, it shows the effectiveness of our algorithm \eqref{alg:constrained} in solving the formulated problem \eqref{prop:main} with exact subgradients.  This observation is consistent with the existing subgradient methods in \cite{gharesifard2014distributed,xi2016distributed,liu2017convergence,zeng2017distributed,zhu2020boundedness}.  
\end{remark}


%
%
%
With Theorem \ref{thm:main-delta}, it is natural for us to enforce some stronger condition on the error $\varepsilon_k$ for a better convergence performance of the entire sequence  $\{x_i(k)\}$.  Along this line, we provide another theorem supposing the accumulated error of subgradient inexactness is not too large.

\begin{theorem}\label{thm:main}
	Suppose Assumptions \ref{ass:convex}--\ref{ass:boundedness} hold. Let the parameters $\alpha_k,\,\e_k>0$ be chosen to satisfy the following condition
	\begin{align}\label{eq:thm-parameter-condition}
		\sum_{k=1}^{\infty}\alpha_k=\infty,\quad \sum_{k=1}^{\infty}\alpha_k^2<\infty,\quad \sum _{k=1}^{\infty} \alpha _{k}\varepsilon _{k}<\infty
	\end{align}
	Then, along the trajectory of algorithm \eqref{alg:constrained}, we have  
	\begin{itemize}
		\item[1)] the sequence $\{||{\bf z}(k+1)-{\bf z}^*||\}$ converges;
		\item[2)] the estimates $x_1(k),\,\ldots,\,x_N(k)$ reach an optimal consensus in the sense that $\lim_{k\to \infty} [{ x}_i(k)-{ x}_j(k)]=0$ and $\lim_{k\to \infty}\tilde f({\bf x}(k))=\tilde f({\bf x}^*)=f^*$;
		\item[3)] $\{{\bf z}(k)\}$  has at least one cluster point  ${\bar{\bf  z}}=\mathrm{col}({\bar {\bf x}},\,{\bf {\bar v}})$ such that $\bar {\bf x}={\bf 1}_N {x^*}$ with ${x^*}$ being  an optimal solution to problem \eqref{prop:main};
		\item[4)] If the optimal solution to problem \eqref{prop:main} is unique, i.e., $\mathcal{X}^*=\{x^*\}$, then $\lim_{k\to \infty} { {x}_i(k)}=x^*$ for each $i\in \N$.
	\end{itemize}
\end{theorem}
\pb  Note that $\Delta({\bf x}(k))\geq 0$ by Lemma \ref{lem:main} and $\sum _{k=1}^{\infty} \alpha _{k}\varepsilon _{k}+\sum_{k=1}^{\infty}\alpha_k^2<\infty$ under the theorem assumption. Applying Lemma 5.31 in \cite{bauschke2017convex} to the inequality \eqref{eq:lem:main}, we can obtain the convergence of $\{||{\bf z}(k+1)-{\bf z}^*||\}$ and 
\begin{align}\label{eq:thm:proof-error-eq5}
	0\leq 	\sum_{i=1}^{_\infty}  \alpha_k\Delta({\bf x}(k))<\infty
\end{align}
Thus, the sequence $\{{\bf z}(k)\}$ must be uniformly bounded by some $C_2>0$. 
From the continuity of $\Delta({\bf x})$, it must be $C_3$-Lipschitz with respect to ${\bf x}$ for some constant $C_3>0$. It follows then 
\begin{align}\label{eq:thm:proof-error-eq6}
	\begin{split}
		\Delta({\bf x}(k+1))-\Delta({\bf x}(k))&\leq C_3||{\bf x}(k+1)-{\bf x}(k)||\\
		&= C_3 ||P_{\tilde X}[{\bf x}(k)- \alpha_k ({\bf g}(k)+L{\bf v}(k)+L{\bf x}(k))]-{\bf x}(k)||\\
		&\leq  C_3||{\bf x}(k)- \alpha_k ({\bf g}(k)+L{\bf v}(k)+L{\bf x}(k))-{\bf x}(k)||\\
		&\leq  C_3  \alpha_k  ||{\bf g}(k)+L{\bf v}(k)+L{\bf x}(k) ||\\
		&\leq  C_3(\sqrt{N}C+2 \lambda_{\max}(L)C_2)\alpha_k
	\end{split}
\end{align}
Jointly using \eqref{eq:thm:proof-error-eq5}, \eqref{eq:thm:proof-error-eq6}, and   $\sum_{k=1}^{ \infty}\alpha_k= \infty$, we resort to Proposition 2 in \cite{alber1998Oprojected} and conclude that $\lim_{k\to \infty} \Delta({\bf x}(k))=0$.  Recalling the expression of $\Delta$,  we have $\lim_{k\to \infty}{\bf x}^\top (k) L{\bf x}(k)=0$ and $\lim_{k\to \infty} [\tilde f({\bf x}(k))+ {{\bf v}^*}^\top L {\bf x}(k)]=\lim_{k\to \infty}\tilde f({\bf x}(k))=\tilde f({\bf x}^*)$. Note that $L$ is positive semidefinite with $0$ as its simple eigenvalue under Assumption \ref{ass:graph}. It follows then $\lim_{k\to \infty} [ x_i(k)-x_j(k)]=0$ and $\lim_{k\to \infty}\tilde f({\bf x}(k))=\tilde f({\bf x}^*)=f^*$.

Due to the uniform boundedness of sequence $\{{\bf z}(k)\}$ by item 1), there must be a convergent subsequence $\{{\bf z}_{k_m}\}$ of $\{{\bf z}_{k}\}$. We denote its limit by $\bar {\bf z}=\mathrm{col}(\bar {\bf x},\,\bar {\bf v})$. Then, it should satisfy $L\bar {\bf x}=0$ and $\tilde {f}(\bar {\bf x})=\lim_{m\to \infty} \tilde f({\bf x}_{k_m})=f^*$ by item 2). In other words, $\bar {\bf x}$ is an optimal solution to \eqref{prop:main-reform}. By Assumption \ref{ass:graph}, one can conclude that there exists some $\bar x\in \R$  such that $\bar {\bf x}={\bf 1}_N \bar x$.  Note that $f(\bar x)=\tilde {f}(\bar {\bf x})=f^*$, i.e., $\bar x$ is an exact optimal solution to problem \eqref{prop:main}.

If $\mathcal{X}^*=\{x^*\}$ holds, all convergent subsequences of $\{{\bf {\bf x}}(k)\}$ have the same limit $x^*$. This combined with the boundedness of $\{{\bf {\bf z}}(k)\}$ implies item 4) and completes the proof. 
\pe

\begin{remark}\label{rem:inexactness}
	This theorem specifies a nontrivial case when our distributed optimization problem \eqref{prop:main} can be exactly solved even using only inexact subgradients information of the local objective functions. This observation is consistent with the centralized results in \cite{kiwiel2004convergence}. 
	Compared with similar primal-dual domain results in \cite{wang2010control,jakovetic2014fast, kia2015distributed, lei2016primal, liu2017constrained, zhu2020boundedness}, this algorithm further allows us to consider nonsmooth objective functions with only approximate subgradients. 
\end{remark}


It is known that normalization might improve the transient performance of the algorithms to avoid overshoots at the starting phase. 
However, conventional normalized techniques often involve some global information and can not be directly implemented in distributed settings.   
We here present a novel componentwise normalized version of algorithm \eqref{alg:constrained} as follows.
\begin{align}\label{alg:constrained-normalized-max}
	x_{i}(k+1) &= P_{X_i}[x_{i}(k)-\frac{\alpha_k}{\max\{c,\,\delta_{ik,D}\}} (g_{i}(k)+\hat x_i(k) + \hat v_i(k))] \nonumber \\
	v_{i}(k+1) & =v_{i}(k)+ \frac{\alpha_k}{\max\{c,\,\delta_{ik,D}\}} \hat x_i(k)\\
	\delta_{ik,\,m}&=\begin{cases} \qquad  ||T^i_{\e_k}({\bf z}(k))||, &\, \mbox{when } m=1 \\ 
		\max_{j\in \mathcal{N}_i}\{\delta_{ik,\,m-1},\, \delta_{jk,\,m-1}\},&\, \mbox{when } 2\leq m\leq D 
	\end{cases}\nonumber
\end{align}
where integer $D \geq \mathrm{D}(\mathcal G)+1$ with $\mathrm{D}(\mathcal{G})$ the diameter of graph $\mathcal{G}$ and $c>0$ is any given constant. Since $\mathrm{D}(\mathcal{G})$ or an upper bound can be computed by distributed rules \cite{oliva2016distributed},  this normalized algorithm is implementable in a fully distributed manner by embedding a max-consensus subiteration. 

Here is a corollary to state that the normalized algorithm \eqref{alg:constrained-normalized-max} retains all the established properties in Theorem \ref{thm:main}.
\begin{corollary}\label{thm:main-normalized}
	Suppose Assumptions \ref{ass:convex}--\ref{ass:graph} hold. Choose the same parameters satisfying condition \eqref{eq:thm-parameter-condition}. Then, along the trajectory of algorithm \eqref{alg:constrained-normalized-max}, the sequence $\{{\bf {\bf z}}(k)\}$ retains all the established properties in Theorem \ref{thm:main}.
	\begin{itemize}
		\item[1)] the sequence $\{||{\bf z}(k+1)-{\bf z}^*||\}$ converges;
		\item[2)] the estimates $x_1(k),\,\dots,\,x_N(k)$ reach an optimal consensus in the sense that $\lim_{k\to \infty} [{\bf x}_i(k)-{\bf x}_j(k)]=0$ and $\lim_{k\to \infty}\tilde f({\bf x}(k))=\tilde f({\bf x}^*)=f^*$;
		\item[3)] $\{{\bf z}(k)\}$ has at least one cluster point ${\bar{\bf  z}}=\mathrm{col}({\bar {\bf x}},\,{\bf {\bar v}})$ such that $\bar {\bf x}={\bf 1}{x^*}$ with ${x^*}$ being  an optimal solution to problem \eqref{prop:main}.
		\item[4)] If the optimal solution to problem \eqref{prop:main} is unique, i.e., $\mathcal{X}^*=\{x^*\}$, then $\lim_{k\to \infty} { {x}_i(k)}=x^*$ for each $i\in \N$.
	\end{itemize}
\end{corollary}
\pb  The proof is similar with that of Theorem \ref{thm:main}. First, we recall Theorem 4.1 in \cite{nejad2009max} and conclude that all agents will get $\max\{c,\,\max_{i\in \N}||T^i_{\e_k}({\bf z}(k))||\}$ after the subiteration. Thus,  we only have to consider the following system:
\begin{align*}
	\begin{split}
		x_{i}(k+1) &= P_{X_i}[x_{i}(k)-\alpha_k\gamma_k (g_{i}(k)+\hat x_i(k) + \hat v_i(k))] \\
		v_{i}(k+1) & =v_{i}(k)+ \alpha_k\gamma_{k} \hat x_i(k)
	\end{split}
\end{align*}
with $\gamma_{k}=\frac{1}{\max\{c,\,\max_{i\in \N}||T^i_{\e_k}({\bf z}(k))||\}}$.  For this new system, we can establish a similar inequality as \eqref{eq:lem:main} for $||{\bf z}(k)-{\bf z}^*||$:
\begin{align*}
	||{\bf z}(k+1)-{\bf z}^*||^2 &\leq ||{\bf z}(k)-{\bf z}^*||^2-2\alpha_k\gamma_{k}({\bf z}(k)-{\bf z}^*)^\top T_{\e_k}({\bf z}(k)) +\alpha_{k}^2\gamma_{k}^2||T_{\e_k}({\bf z}(k))||^2
\end{align*}
Note that $\gamma_{k}^2||T_{\e_k}({\bf z}(k))||^2\leq N$ and $0\leq \gamma_{k}\leq \frac{1}{c}$. Recalling the fact \eqref{eq:thm-proof-error-eq2} and $\Delta({\bf x}(k))\geq 0$, one can obtain
\begin{align*}
	||{\bf z}(k+1)-{\bf z}^*||^2&\leq ||{\bf z}(k)-{\bf z}^*||^2-2\alpha_k\gamma_{k}\Delta({\bf x}(k)) +\frac{2N}{c} \alpha_{k} \e_k+N\alpha_{k}^2
\end{align*}
According to  Lemma 5.3.1 in \cite{bauschke2017convex}, we can conclude the convergence of $\{||{\bf z}(k+1)-{\bf z}^*||\}$ and $\sum_{k=1}^{\infty} \alpha_k\gamma_{k}\Delta({\bf x}(k))\leq \infty $. Then, $\{||{\bf z}(k)||\}$ is uniformly bounded, which implies that there must be small enough constant $C_4>0$ such that $\gamma_{k}\geq C_4>0$.  We use Proposition 2 in \cite{alber1998Oprojected} again and obtain that 
$\lim_{k\to \infty} \Delta({\bf x}(k))=0$.
Then, items 3) and 4) can be easily verified following a similar procedure as in Theorem \ref{thm:main}. The proof is thus complete.
\pe

\begin{remark}
	Compared with the conventional normalized step sizes  in \cite{ruszczynski2011nonlinear, boyd2014subgradient}, the proposed componentwise normalized step size can be taken as their distributed extension and the iterative sequence generated by \eqref{alg:constrained-normalized-max}  might have a better transient behavior than that generated by \eqref{alg:constrained}.  Interestingly, the widely used subgradient boundedness assumption (i.e., Assumption \ref{ass:boundedness})  is also removed as a byproduct, which might be favorable in distributed scenarios.
\end{remark}

\section{Simulation}\label{sec:simu}

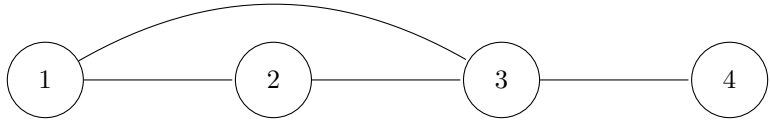
\begin{figure}
	\centering
	\begin{tikzpicture}[shorten >=1pt, node distance=3.0 cm, >=stealth',
		every state/.style ={circle, minimum width=  1 cm, minimum height= 1 cm}]
		\node[align=center,state](node1) {   1};
		\node[align=center,state](node2)[right of=node1]{   2};
		\node[align=center,state](node3)[right of=node2]{  3};
		\node[align=center,state](node4)[right of=node3]{   4};
		\path[-]  (node1) edge (node2)
		(node2) edge (node3)
		(node3) edge (node4)
		(node1) edge [bend left=30]  (node3)
		;
	\end{tikzpicture}
	\caption{Communication graph $\mathcal G$ in our example.}\label{fig:graph}
\end{figure}

In this section, we consider a LASSO (least absolute shrinkage and selection operator) regression problem to verify the effectiveness of our algorithms:
\begin{align*}
	\min_{x\in X} f(x)=\frac{1}{2}\sum_{i=1}^N ||x-p_i||^2+\lambda N ||x||_1
\end{align*}
where $\lambda N>0$ is the regularization parameter, $p_i$ is an estimate only known by agent $i$, and $X$ is the constrained set. In this case, we let $X_i=X$ and $f_i(x)=\frac{1}{2}||x-p_i||^2+\lambda ||x||_1$ and put it into a form of  problem \eqref{prop:main}.  
Distributed subgradient  algorithms have been developed to solve this problem when $\partial f_i$ is available, e.g., \cite{liu2017convergence}. Although $\partial f_i$  and  $\partial_{\varepsilon} f_i$ can be easily calculated, we use this example to show the effectiveness of our algorithm only using its $\varepsilon$-subdifferential instead of the exact one.

According to the definition of $\varepsilon$-subgradient, we have
$$ 
\partial_\varepsilon f_i(x)=
\begin{cases}
	[x-p_i-\lambda,x-p_i-\lambda-\frac{\lambda \varepsilon}{x}] & \text{for   }  x<-\frac{\varepsilon}{2}, \\
	[x-p_i-\lambda,x-p_i+\lambda]  & \text{for   }  x \in [-\frac{\varepsilon}{2},\frac{\varepsilon}{2}] ,\\
	[x-p_i+\lambda-\frac{\lambda \varepsilon}{x},x-p_i+\lambda] &\text{for   } x>\frac{\varepsilon}{2}
\end{cases}   
$$

\begin{figure}
	\centering
	\subfigure{
		\includegraphics[width=0.42\textwidth]{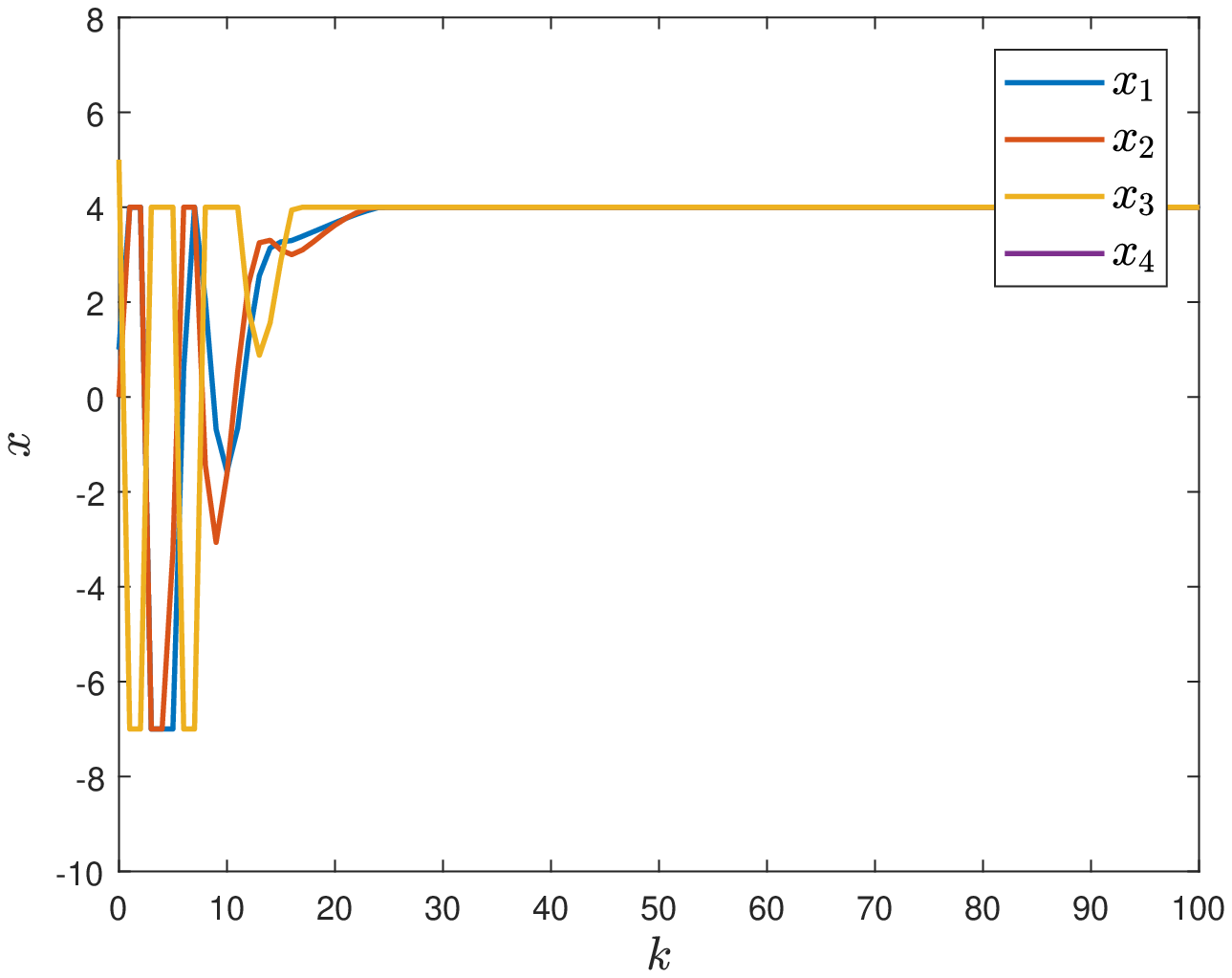}
	} 
	\subfigure{
		\includegraphics[width=0.42\textwidth]{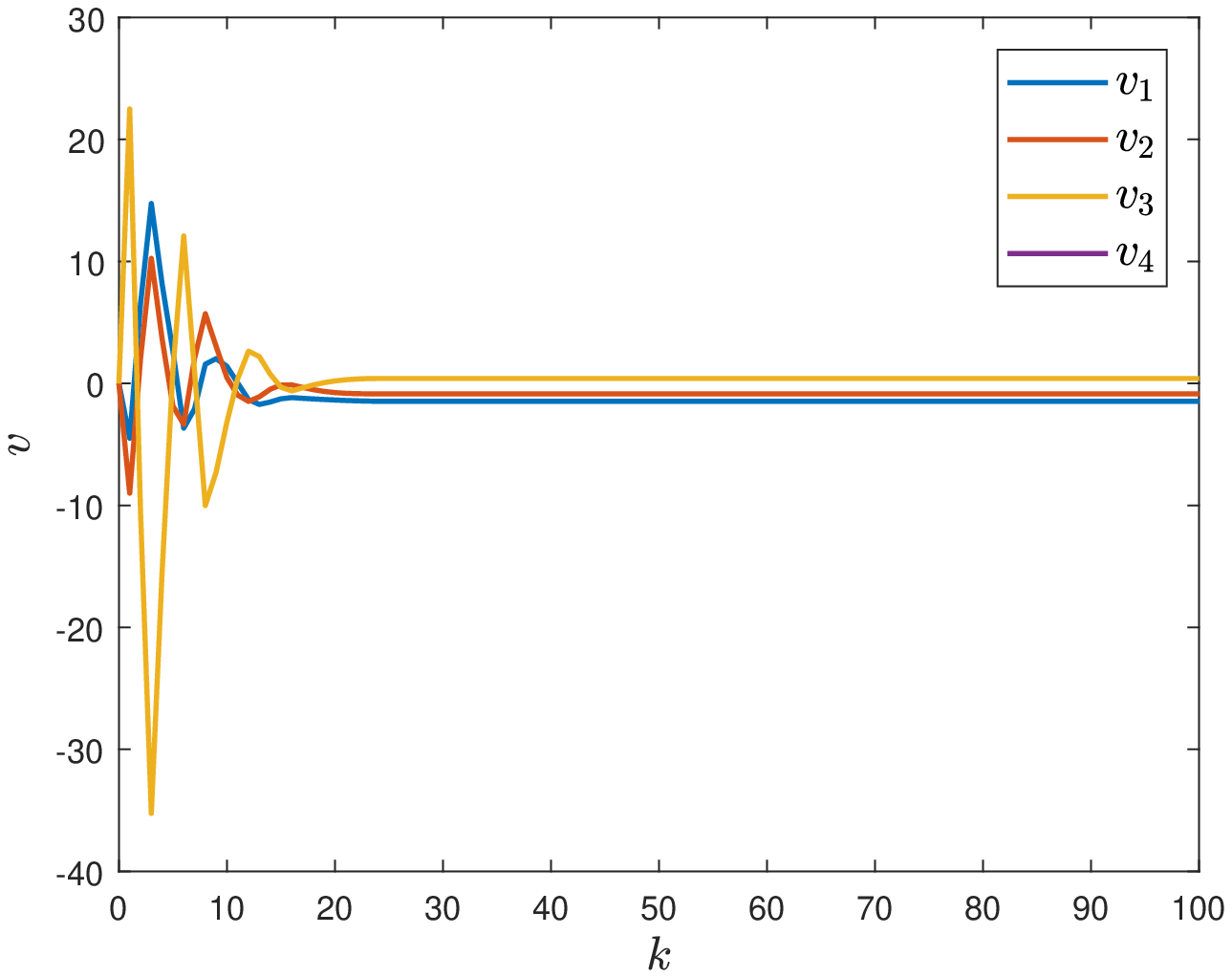}
	}
	\caption{Profiles of primal and dual variables in Algorithm \eqref{alg:constrained}.}\label{fig:simu-1}
\end{figure}

\begin{figure}
	\centering
	\subfigure{
		\includegraphics[width=0.42\textwidth]{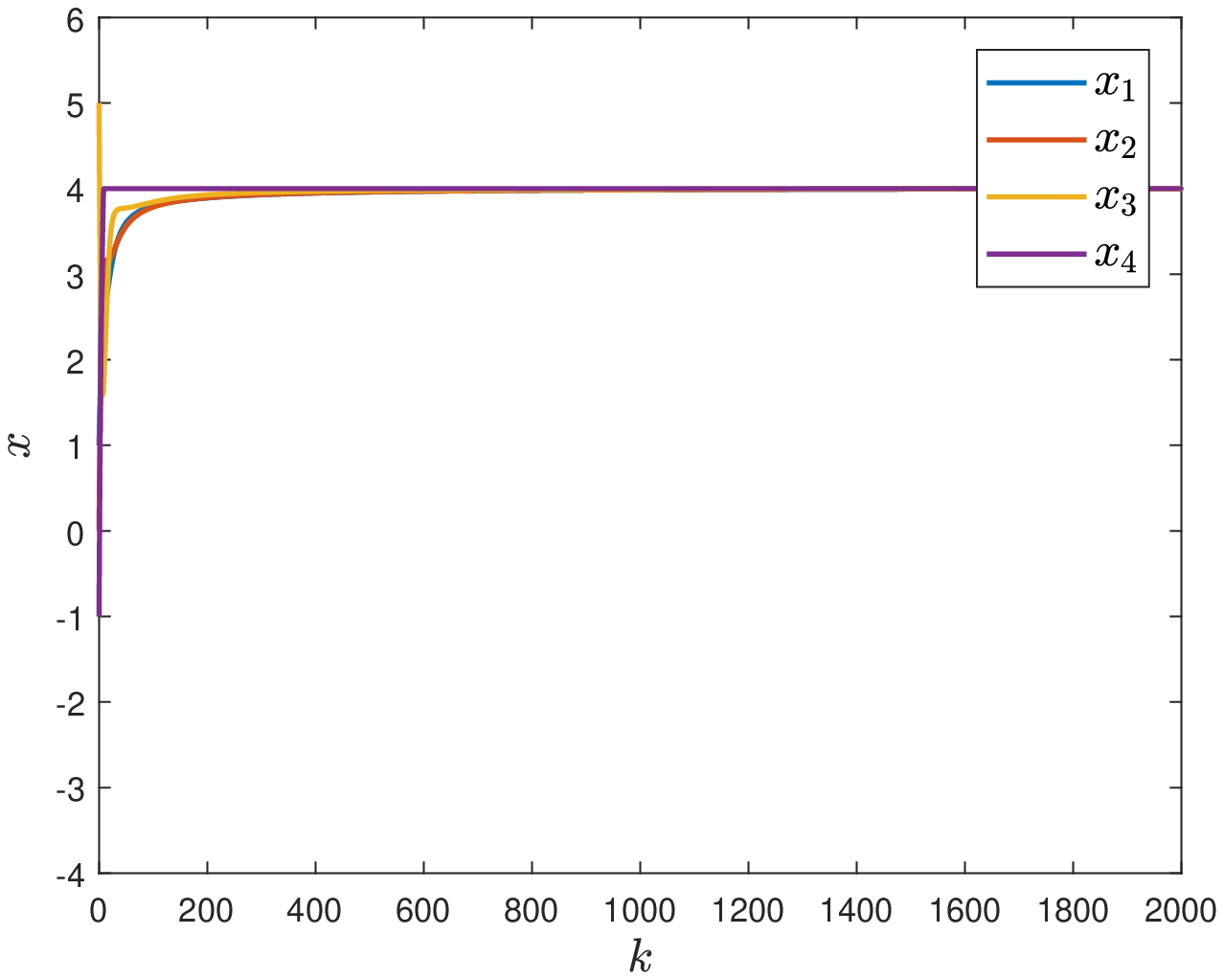}
	} 
	\subfigure{
		\includegraphics[width=0.42\textwidth]{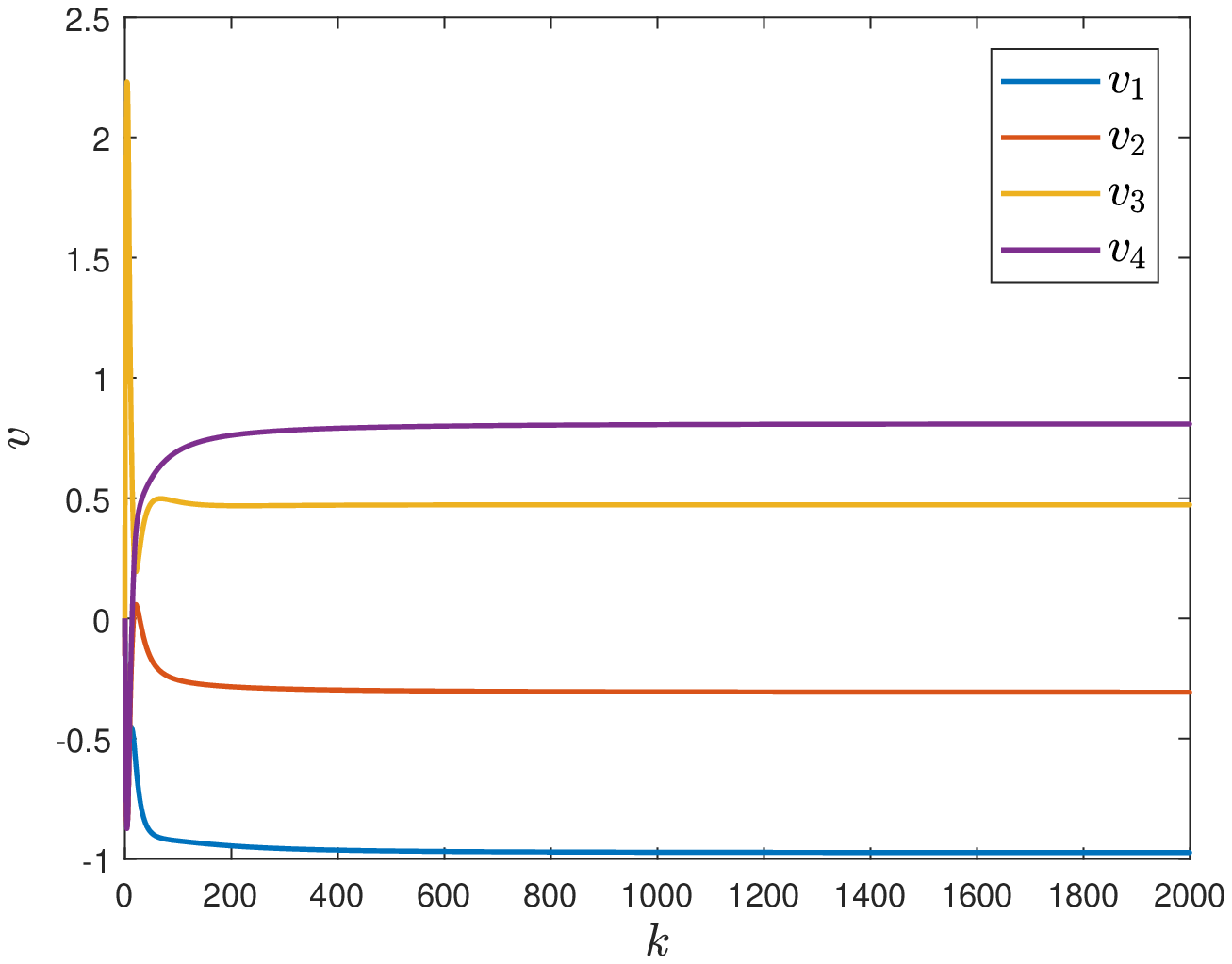}
	}
	\caption{Profiles of primal and dual variables in Algorithm \eqref{alg:constrained-normalized-max}.}\label{fig:simu-2}
\end{figure}

For simulations, we let $N=4$, $p_i=2i$, and $\lambda=0.1$.  The communication graph  is given as Fig.~\ref{fig:graph} with unity weights. To make it more interesting, we assume this problem has set constraints specified by $X_i=[-11+i,\,8-i]$ with $i=1,\,2,\,3,\,4$. Assumptions \ref{ass:convex}-\ref{ass:graph} can be confirmed. 
To ensure the condition \eqref{eq:thm-parameter-condition}, we choose $\alpha_k=\varepsilon_k=\frac{3}{k+1}$.  Then, the considered problem can be solved by our proposed distributed primal-dual $\varepsilon$-subgradient method (PD$\varepsilon$SM) \eqref{alg:constrained} and the normalized primal-dual $\varepsilon$-subgradient method (NPD$\varepsilon$SM) \eqref{alg:constrained-normalized-max} according to Theorem  \ref{thm:main} and Corollary \ref{thm:main-normalized}.   In simulations, when $x<-\frac{\varepsilon}{2}$, we choose $x-p_i-\lambda-\frac{\lambda \varepsilon}{x}$ as the $\varepsilon$-subgradient, when $-\frac{\varepsilon}{2}\leq x\leq \frac{\varepsilon}{2}$, we choose $x-p_i+\lambda$ as the $\varepsilon$-subgradient, when $x>\frac{\varepsilon}{2}$, we choose $x-p_i+\lambda-\frac{\lambda \varepsilon}{x}$ as the $\varepsilon$-subgradient.

\begin{figure}
	\centering
	\includegraphics[width=0.84\textwidth]{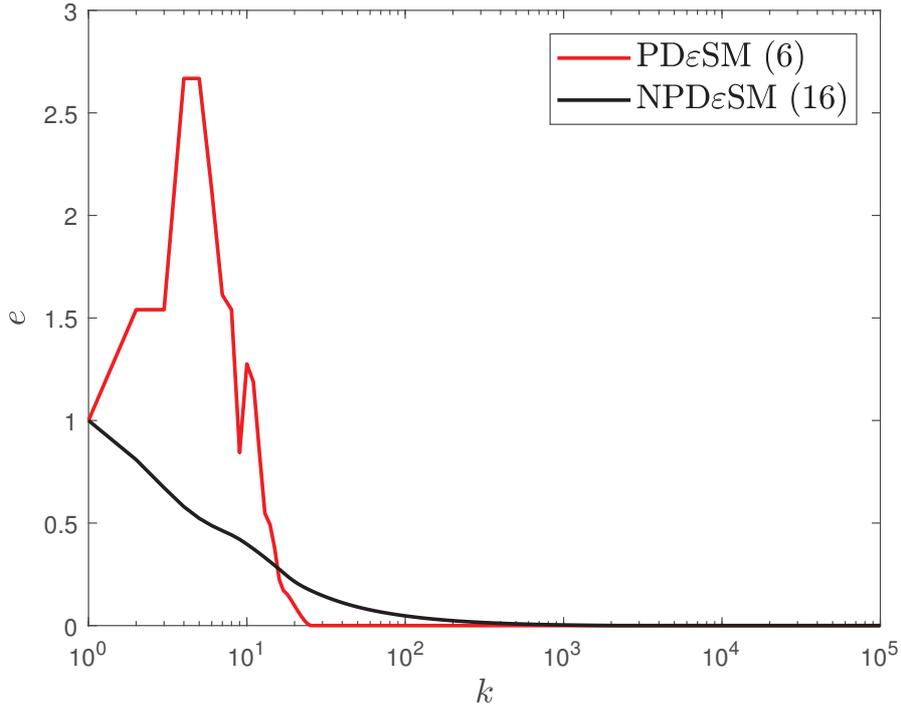}
	\caption{Profiles of residential errors in Algorithms \eqref{alg:constrained} and \eqref{alg:constrained-normalized-max}.}\label{fig:simu-3}
\end{figure}

Simulation  results  with ${\bf x}(1)=\mathrm{col}(1,\, 0,\, 5,\, -1)$ and $c=0.1$ for \eqref{alg:constrained} and \eqref{alg:constrained-normalized-max} are shown in Figs.~\ref{fig:simu-1}--\ref{fig:simu-2}, where all agents' primal variables are observed to converge to the global optimal solution $x^*=4$ while the dual variables are bounded and converge. This verifies the effectiveness of our algorithms. Moreover, one can find that although the normalized step size in \eqref{alg:constrained-normalized-max} might slow down the convergence speed compared with \eqref{alg:constrained}, the resultant transient performance of the primal and dual variables has been much improved with less and weaker oscillations.  For a clear comparison, we let $e(t)=\frac{||x_i(k)-{\bf 1}_4 x^*||}{||x_i(1)-{\bf 1}_4 x^*||}$ be the residential error of our algorithms. The profiles of $e(k)$ in both algorithms are shown in  Fig.~\ref{fig:simu-3}. From this, we can also confirm the improvement of transient performance by the proposed componentwise normalized step size.

\section{Conclusion}\label{sec:con}
In this paper, we have attempted to solve a distributed constrained optimization problem with inexact subgradient information of local objective functions. We have developed a projected primal-dual dynamics using only $\varepsilon$-subgradients and discussed its convergence properties. In particular, we have shown the exact solvability of this problem if the accumulated error introduced by  subgradient inexactness is not too large.  We have also presented a novel distributed normalized step size to improve the transient performance of our algorithms. It is interesting to consider more general graphs in future works. 

\vspace{6pt}




\end{document}